\documentclass{amsproc}

 \newtheoremstyle{plain}
     {3pt}
     {3pt}
     {}
     {}
     {\sc}
     {:}
     {.5em}
     {}

 \newtheoremstyle{tshtheorem}
     {3pt}
     {5pt}
     {\it}
     {}
     {\sc}
     {.}
     {.5em}
     {}

\theoremstyle{tshtheorem}
\newtheorem{theorem}{Theorem}[section]

\newtheorem{thmd}[theorem]{Darboux's Theorem}
\newtheorem{thmc}[theorem]{Convexity Theorem}
\newtheorem{thmi}[theorem]{Injectivity Theorem}
\newtheorem{thms}[theorem]{Surjectivity Theorem}
\newtheorem{thmcr}[theorem]{CR Surjectivity Theorem}

\newtheorem{definition}[theorem]{Definition}

\theoremstyle{plain}
\newtheorem{example}[theorem]{Example}

\newtheorem{remark}[theorem]{Remark}

\numberwithin{equation}{section}
\numberwithin{figure}{section}

\usepackage{epsfig}
\usepackage{psfrag} 
\input xy
\xyoption{all}

\newcommand{\lcm}{\mathrm{lcm}}

\newcommand{\C}{{\mathbb{C}}}

\newcommand{\Z}{{\mathbb{Z}}}
\newcommand{\Q}{{\mathbb{Q}}}
\newcommand{\R}{{\mathbb{R}}}

\newcommand{\X}{{\mathfrak{X}}}
\newcommand{\XX}{{\mathcal{X}}}
\newcommand{\Y}{{\mathcal{Y}}}
\newcommand{\ZZ}{{\mathcal{Z}}}

\newcommand{\algt}{\mathfrak{t}}
\newcommand{\algg}{\mathfrak{g}}

\newcommand{\iso}{\cong}

\newcommand{\Year}{{a}}
\newcommand{\im}{\mathrm{im}}
\newcommand\AGV{{$\aleph$}\!\! GV}
\newcommand{\age}{\mathrm{age}}

\begin{document}

\title[Cohomology of abelian symplectic reductions]{Orbifold cohomology of abelian symplectic reductions and the case of weighted projective spaces}
\author{Tara S. Holm}
\address{Department of Mathematics, Cornell University, Ithaca, NY 14853-4201 USA}
\email{tsh@math.cornell.edu}
\thanks{TSH is grateful for the support of the NSF through the grant DMS-0604807.}

\subjclass{Primary 53D20; Secondary 14N35, 53D45, 57R91}

\keywords{Symplectic quotient, orbifold, cohomology}

\begin{abstract}
These notes accompany a lecture about the topology of symplectic (and
other) quotients.  The aim is two-fold: first to advertise
the ease of computation in the symplectic category; and second to 
give an account of some new computations for weighted projective spaces.
We start with a brief exposition of how orbifolds arise in the 
symplectic category, and discuss the techniques used to understand their
topology. We then show how these results can be used to compute the 
Chen-Ruan orbifold cohomology ring of
abelian symplectic reductions.  We conclude by comparing the several rings
associated to a weighted projective space.  We make these computations
directly, avoiding any mention of a stacky fan or of a labeled moment polytope.
\end{abstract}

\maketitle

\tableofcontents

\renewcommand{\arraystretch}{1.6}

The notion of an {\bf orbifold} has been present in topology since the 1950's \cite{Sat56,Sat57}.   More recently, orbifolds have played an important role in differential and algebraic geometry, and in mathematical physics.  A fundamental theme is to compute topological invariants associated to an orbifold, with one ostensible goal to understand Gromov-Witten invariants for these spaces.  The aim of the present article is modest: to expound how techniques from symplectic geometry may be used to understand the degree-zero genus-zero Gromov-Witten invariants with three marked points, the so-called {\bf Chen-Ruan orbifold cohomology ring}; and to make explicit the details of these techniques in the case of weighted projective spaces.
 
In the symplectic category, orbifolds
arise as symplectic quotients. We recount the techniques from
symplectic geometry that may be used to compute topological
invariants of a symplectic quotient.  This is based on Kirwan's
seminal work \cite{K:quotients}; and for orbifold invariants, the author's joint work
with Goldin and Knutson \cite{GHK:preorb}.   The quotients 
we consider are by a compact connected abelian group.  We employ techniques coming from algebraic topology, most notably using equivariant cohomology.  For those used to working with finite groups, it is important to note that, whereas for finite groups the 
invariant part of a cohomology ring is identical to the equivariant cohomology, this is 
not the case for connected groups.

The main example in this article is 
a {\bf weighted projective space} $\C P^n_{(b)}$.  Its definition depends on a 
sequence $(b) = (b_0,\dots,b_n)$ of positive integers.  Kawasaki showed that 
the ordinary cohomology groups, with integer coefficients, of
the underlying topological space of a weighted projective space  are identical to the 
cohomology groups of a 
smooth projective space 
\cite{K:weightedproj}, but there is a twisted ring structure.  We review the details 
of his work.  Then in 
Theorem~\ref{th:surjZwps}, we compute the cohomology of the orbifold $[\C P^n_{(b)}]$,
proving that
\begin{equation}\label{eq:wps}
H^*([\C P^n_{(b)}];\Z) = H_{S^1_{(b)}}^*(S^{2n+1};\Z) \cong \frac{\Z[u]}{\langle b_0\cdots b_n u^{n+1}\rangle}.
\end{equation}
Whereas Kawasaki finds a twist in the the ring structure, we find torsion in high degrees of the ring \eqref{eq:wps}.
There is a natural map from Kawasaki's ring to this one, and we describe the map explicitly.
Finally in Theorem~\ref{th:cr}, we compute the Chen-Ruan cohomology ring of this orbifold.
We make this computation using integer coefficients, generalizing results in \cite{J:wps,Mann:thesis, Mann:paper}.  Moreover, we give explicit generators and relations, and avoid  mentioning a stacky fan \cite{BCS:toricvarieties} or a labeled polytope \cite{LT:toricorbifolds,GHK:preorb}.

The definitions in this article make sense for arbitrary coefficient rings.  Indeed, all computations in the final section use integer coefficients.  Moore and Witten have suggested that the torsion in $K$-theory has more physical significance than torsion in cohomology \cite{MoWi:torsion}.  The  author together with Goldin, Harada and Kimura, is investigating a $K$-theoretic version of \cite{GHK:preorb} and of the computations herein, building on the work of Harada and Landweber \cite{HL}.

The remainder of the paper is organized as follows.  In Section~\ref{sec:symp1} we give a quick exposition of how orbifolds arise in the symplectic category.  We then introduce several cohomology rings associated to an orbifold in Section~\ref{sec:orb}.  We advertise the ease of computation for these rings in Section~\ref{sec:symp}.  The novel results in this article are the computations in Section~\ref{sec:wps}.  We include detailed proofs that avoid much of the symplectic machinery used in \cite{GHK:preorb}.

 \medskip
 
 \noindent {\bf Acknowledgments.}
Many thanks are due to Tony Bahri, Matthias Franz, Rebecca Goldin, Megumi Harada,
Ralph Kaufmann, Takashi Kimura, Allen Knutson, Eugene Lerman, Reyer Sjamaar, and Alan Weinstein for many helpful conversations; and to  Yoshiaki Maeda and the organizers and sponsors of Poisson 2006 in Tokyo, Japan, where this work was presented.

\section{Symplectic manifolds and quotients}\label{sec:symp1}

We begin with a very brief introduction to the symplectic category; a
more detailed account of the subject can be found in \cite{CdS:book}.
A {\bf symplectic form}
on a manifold $M$ is a closed non-degenerate two-form $\omega\in
\Omega^2(M)$.  Thus, for any tangent vectors $\XX,\Y\in T_pM$,
$\omega_p(\XX,\Y)\in\R$.  The key examples include the following.

\begin{example}\label{eg:2sphere}
$M=S^2 = \C P^1$ with $\omega_p (\XX,\Y)$ equal to the signed area of the
parallelogram spanned by $\XX$ and $\Y$.  This is the Fubini-Study form
on $\C P^1$. 
\end{example}

\begin{example}
$M$ any orientable Riemann surface with $\omega$ as in Example~\ref{eg:2sphere}.  Note that orientability is a necessary condition on a symplectic manifold $M$, because the top exterior power of the symplectic form is a volume form.
\end{example}

\begin{example}\label{eg:darboux}
$M= \R^{2d}$ with $\omega = \sum dx_i\wedge dy_i$.
\end{example}

\begin{example}
$M= \mathcal{O}_\lambda$ a coadjoint orbit of a compact
connected semisimple Lie group, equipped with $\omega$ the
Kostant-Kirillov-Soriau form. 
\end{example}

\noindent Example~\ref{eg:darboux} gains particular importance because of 
\begin{thmd}
Let $M$ be a symplectic $2d$-manifold with symplectic form $\omega$.
Then for every point $p\in M$, there exists a coordinate chart $U$
about $p$ with coordinates $x_1,\dots,x_d,y_1,\dots,y_d$ so that on
this chart, 
\begin{equation}
\label{ }
\omega = \sum_{i=1}^d dx_i\wedge dy_i.
\end{equation}
\end{thmd}
\noindent Thus, whereas Riemannian geometry uses local invariants
such as {\bf curvature} to distinguish metrics, symplectic forms are locally
indistinguishable. 

The symmetries of a symplectic manifold may be encoded as a group action.  Here we restrict ourselves to a compact connected abelian group $T= (S^1)^n$.  An action of $T$ on $M$ is {\bf symplectic} if it preserves $\omega$; that is, $\rho_g^*\omega = \omega$, for each $g\in T$, where $\rho_g$ is the diffeomorphism corresponding to the group element $g$.  The action is {\bf Hamiltonian} if in addition, for every $\xi\in\algt$, the vector field
\begin{equation}
\label{ }
\XX_\xi = \frac{d}{dt}[ \exp (t\xi)] |_{t=0}
\end{equation}
is a Hamiltonian vector field.  That is, we require that 
$
\omega(\XX_\xi, \cdot ) = d\phi^\xi
$
is an exact one-form.  Each $\phi^\xi$ is a smooth function on $M$, determined up to a constant.  Taking them together, we may define a {\bf moment map}
\begin{equation}
\label{ }
\begin{array}{rcc}
\Phi: M  & \longrightarrow & \algt^* \\
 p & \longmapsto & \left(\begin{array}{rcl}
 				\Phi(p):\algt & \rightarrow & \R \\
				\xi & \mapsto & \phi^\xi(p)
				\end{array}\right).
\end{array}
\end{equation}
Returning to our examples, we have Hamiltonian actions in all but the
second example.

\begin{example}
The circle $S^1$ acts on $M=S^2 = \C P^1$ by rotations.  If we use
angle and height coordinates on $S^2$, then the vector field this
action generates is tangent to the latitude lines, so in coordinates, 
$\XX^\xi = \frac{\partial}{\partial\theta}$, and since $\omega =
d\theta\wedge dh$, $\omega(\XX_\xi, \cdot ) = dh$, so a moment map is
the height function on $S^2$, as shown in Figure~\ref{fig:sphere} below.  
\begin{figure}[h]
  \begin{center}
  \psfrag{S}{$S^1$}
  \psfrag{Phi}{\small $\Phi$}
    \epsfig{file=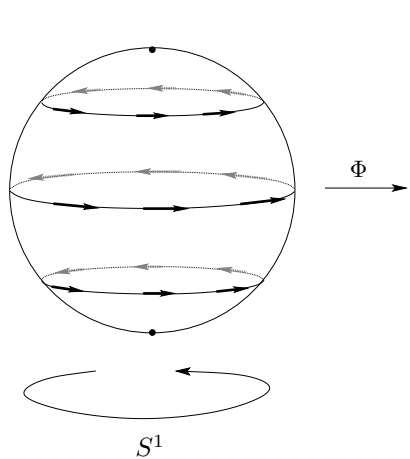,height=6cm}
    \caption{The vector field and moment map for $S^1$ acting by
  rotations on $S^2$.} \label{fig:sphere}
  \end{center}
 \end{figure}
\end{example}

\begin{example}
If $M$ is a two-torus $M=T^2 = S^1\times S^1$, then $S^1\times
S^1$ acts on itself by multiplication.  This action is symplectic, but
is not Hamiltonian.  In fact, no Riemann surface with non-zero genus has a
nontrivial Hamiltonian torus action. 
\end{example}

\begin{example}
The torus $T^d=(S^1)^d\subset \C^d$ acts by coordinate-wise multiplication on $M= \R^{2d} = \C^d$. 
This action rotates each copy of $\C = \R^2$ (at unit speed), and is Hamiltonian. Identifying $\algt^*
\cong \R^d$, a moment map is 
\begin{equation}
\label{ }
\Phi (z_1,\dots,z_n) = (|z_1|^2, \dots,|z_d|^2),
\end{equation}
up to a constant multiple.
\end{example}

\begin{example}
Each coadjoint orbit $M= \mathcal{O}_\lambda\subseteq \algg^*$
may be identified as a homogenous space $G/L$, where $L$ is a Levi
subgroup of the Lie group $G$.  Thus $G$ and its maximal torus
$T$ act on $M$ by left multiplication.  A $G$-moment map is
inclusion 
\begin{equation}
\label{ }
\Phi_G : \mathcal{O}_\lambda\hookrightarrow \algg^*,
\end{equation}
and a $T$-moment map is the $G$-moment map composed with the natural
projection $\algg^* \to\algt^*$ that is dual to the inclusion
$\algt\hookrightarrow\algg$.  
\end{example}

In each of these examples, the image of the (torus) moment map is a convex
subset of $\R^n$.  This is true more generally.
\begin{thmc}[\cite{At:convexity},\cite{GS:convexity}]\label{thm:convex} 
If $M$ is a compact Hamiltonian $T$-space, then $\Phi(M)$ is a convex
polytope.  It is the convex hull of $\Phi(M^T)$, the images of the
$T$-fixed points.
\end{thmc}
The convexity theorem is an example of a {\bf localization
phenomenon}: a global feature (the image of the moment map) that is
determined by local features of the fixed points (their images under
the moment map).  The convexity property  is a recurring theme in
symplectic geometry; its many guises are illustrated in
\cite{GuSj:convex}.

The moment map is a $T$-invariant map: it maps entire $T$-orbits to
the same point in $\algt^*$.  Thus when $\alpha$ is a regular value, the level set
$\Phi^{-1}(\alpha)$ is a $T$-invariant submanifold of $M$.  Moreover, the action of 
$T$ on a regular level set
is {\bf locally free}: it has only finite stabilizers.  This follows
directly from the moment map condition: at a regular value,
$d\phi^\xi$ is never zero, implying that $\XX_\xi$ is not zero, so
there is no $1$-parameter subgroup fixing points in the level set.
Thus, at a regular value the {\bf symplectic reduction} $M/\!\!/T(\alpha) =
\Phi^{-1}(\alpha)$ is an orbifold.  In fact, Marsden and Weinstein proved
\begin{theorem}[\cite{MW:reduction}]
If $M$ is a Hamiltonian $T$-space and $\alpha$ is a regular value of the 
moment map $\Phi$, then the symplectic reduction $M/\!\!/T(\alpha)$ is 
a {\bf symplectic} orbifold. 
\end{theorem}

\noindent More generally, the symplectic reduction $M/\!\!/T(\alpha)$ at
a critical value is a symplectic stratified space \cite{SL:strat}. 

Symplectic reduction is an important technique for constructing new
symplectic manifolds from old. From our examples, we may construct
several classes of symplectic manifolds.

\begin{example}
For the action of $S^1$ on $M=S^2 = \C P^1$ by rotation, the
level set of a 
regular value is a latitude line, which the circle rotates.  The
quotient is a point.  Note that if $S^1$ acts by rotation $S^2$ at
twice the usual speed, then the quotient, as an orbifold, is $[
\mathrm{pt}/\Z_2]$.  Thus, every  orbifold $[ \mathrm{pt}/\Z_k]$ is a
quotient of $S^2$ by a rotation action.
\end{example}

\begin{example}
$T^d$ acts on $M= \R^{2d} = \C^n$ by rotation of each copy of
$\C = \R^2$.  The level set of a regular value is a copy of $T^d$, and
so again the quotient is a point (or potentially an orbipoint, for
different actions of $T^d$).  However, we may also restrict our
attention to subtori $K\subseteq T^n$.  The action of $K$ is still
Hamiltonian, and for certain choices of $K$, $\C^n/\!\!/K (\alpha)$ is a
{\bf symplectic toric orbifold}.  Lerman and Tolman show that every
{\bf effective}  symplectic toric orbifold may be constructed in this way
\cite{LT:toricorbifolds}.  
\end{example}

\begin{example}
For the $T$-action on a coadjoint orbit $M=
\mathcal{O}_\lambda$, the symplectic reduction $M/\!\!/T(\alpha)$ is
known as a {\bf weight variety}.  One may determine the possible
orbifold singularities by analyzing the combinatorics of $G$ and its
Weyl group.  See, for instance, \cite{AllenThesis} or  \cite{GHK:preorb}.
This reduced space plays an important role in representation theory.
\end{example}

\section{Orbifolds and their cohomology}\label{sec:orb}

We now turn to orbifolds in the topological category.  In terms of
local models, an {\bf orbifold} is a topological space where each
point has a neighborhood homeomorphic (or diffeomorphic) to the
quotient of a (fixed dimensional) vector space by a finite
group. Satake introduced  this notion in the 1950's
\cite{Sat56,Sat57}, originally calling the spaces $V$-manifolds.
Thurston coined the term orbifold when he rediscovered them in the
1970's (see \cite{Thu}) in his study of $3$-manifolds.  This local model, however, makes it
difficult to define very basic pieces in the theory of orbifolds:
overlap conditions on orbifold charts, suborbifolds, and maps between
orbifolds.  As is evident already in the work of Haefliger
\cite{haefliger}, the proper way to think of an orbifold is as a
Morita equivalence class of {\bf groupoids}, one of which is a proper
\'etale groupoid (see \cite{Moe}); or equivalently as a smooth {\bf Deligne-Mumford
stack} (see \cite{DM}). For example, using this structure, a map of orbifolds should
simply be a morphism of the appropriate objects. 

While groupoids or stacks provide the correct mathematical framework,
the technology is a bit beyond the scope of this article.  Indeed, for
us it is sufficient to work with the local models, largely
because we restrict our attention to orbifolds that arise as global
quotients.  Nevertheless, we will need to distinguish between an orbifold $\X$ or $[X]$ and its
underlying topological space (or {\bf coarse moduli space}) $X$. In particular, when $\X$ is
presented as a global quotient of a manifold $M$ by a group $G$, 
we will use square brackets $[M/G]$ to denote the
orbifold, and $M/G$ to denote the underlying coarse moduli.

For an orbifold 
$\X$, at each  point $x\in \X$, we
have a local isotropy group $\Gamma_x$ at $x$.  We will be interested
in {\bf almost complex orbifolds}, that is orbifolds that have local
models isomorphic to $\C^d/ \Gamma_x$, with $\Gamma_x\subseteq U(d)$ a
finite group acting  unitarily on $\C^d$. Our main example in this
section is the orbisphere shown in the figure below.   
\begin{figure}[h]
  \begin{center}
    \epsfig{file=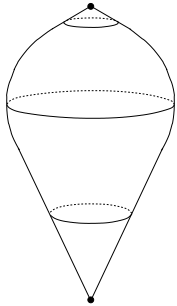,height=4cm}
    \caption{An orbisphere with two orbifold points, a $\Z_p$
    singularity at the north pole and a $\Z_q$ singularity at the
    south pole.  All other points in the
    space have local isotropy group the one-element group.
    When $p$ and $q$ are relatively prime, this is a
    weighted projective space $\C P_{p,q}^1$.  
    }\label{fig:orbisphere}  	
  \end{center}
 \end{figure}
This is a
symplectic toric orbifold in the sense of Tolman and Weitsman
\cite{LT:toricorbifolds}.  In that context, it corresponds to the
labeled polytope that is an edge, with the vertices labeled $p$ and
$q$. 
This orbifold cannot (always) be presented as a global quotient by a
{\bf finite group}, although it can be presented as a symplectic
reduction.  When $p$ and $q$ are relatively prime, it is the reduction
of $\C^2$ by the Hamiltonian $S^1$-action 
\begin{equation}
\label{ }
e^{2\pi i\theta}\cdot (z_1,z_2) = (e^{p\cdot 2\pi i\theta}\cdot z_1,
e^{q\cdot 2\pi i\theta}\cdot z_2).
\end{equation}
Thus, it is a global quotient of the level set
$\Phi^{-1}(\alpha)\approx S^3$ by a locally free $S^1$ action. In this
case, it is also the global quotient of $\C P^2$ by a $\Z_p\times
\Z_q$ action\footnote{With an apology to number theorists,
 we take the topologist's notation:  $\Z_p$ denotes the integers
modulo $p$.}.  When $p$ and $q$ are not relatively prime, the
orbisphere is not a global quotient by a finite group, but it is still
a symplectic quotient of $\C^2$ by a Hamiltonian $S^1\times \Z_g$
action, where $g=\gcd(p,q)$ (following \cite{LT:toricorbifolds}).   Note also that an orbisphere is 
isomorphic to a weighted projective space $\C P^2_{p,q}$ exactly when
$p$ and $q$ are relatively prime.  If $p$ and $q$ are not relatively prime,
the weighted projective space is not {\bf reduced}: it has a global stabilizer.
On the other hand,
the orbisphere described above is always reduced.

Next we turn to algebraic invariants that we may attach to an orbifold
$\X$.  The hope is that these invariants are computable and at the same
time retain some information of the orbifold structure of $\X$.  All
the invariants are isomorphic to singular cohomology if $\X=X$ is in fact
a manifold. 
\begin{definition}
The {\bf ordinary cohomology ring} of an orbifold $\X$ is the singular
cohomology of the underlying topological space $X$, 
\begin{equation}
\label{ }
H^*(X;R),
\end{equation}
with coefficients in a commutative ring $R$.  
\end{definition}

\noindent This ring is computable using standard techniques from
algebraic topology, but it does not distinguish between the orbisphere
in Figure~\ref{fig:orbisphere} and a smooth sphere.   

For the second invariant, we restrict our attention to orbifolds $\X$
presented as the quotient of a manifold $M$ by the locally free action
of a Lie group $G$.  It is conjectured that every orbifold can be
expressed as such a global quotient, and it is known to be true for
{\bf effective} or {\bf reduced orbifolds}, those that do not have a
global finite stabilizer (see, for example, \cite[Theorem
2.18]{EHKV}).  A presentation as a global quotient is desirable
because then the topology or geometry of the quotient $\X$ is simply
the $G$-{\bf equivariant} topology or geometry of the manifold
$M$. This principle motivates the following definition.

\begin{definition}
Given a presentation of an orbifold $\X = [M/G]$ as a global quotient,
the {\bf cohomology ring of the orbifold} $\X$ is the equivariant
cohomology ring 
\begin{equation}
\label{ }
H^*(\X;R) := H^*_G(M;R),
\end{equation}
with coefficients in a commutative ring $R$.  
\end{definition}
\noindent Recall that {\bf equivariant cohomology} is a generalized
cohomology theory in the equivariant category.  Using the {\bf Borel
model}, we define 
\begin{equation}
\label{ }
H_G^*(M;R) := H^*((M\times EG)/G;R),
\end{equation}
where $EG$ is
a contractible (though infinite dimensional) space with a free $G$
action, and $G$ acts diagonally on $M\times EG$.  There are
well-developed methods for computing equivariant cohomology, hence
this invariant is still computable.

Whenever $\X=[M/G]$ is a global quotient, the associated quotient map 
\begin{equation}
q:M\to X
\end{equation}
is a $G$-invariant map. This induces a continuous map 
\begin{equation}
q: (M\times EG)/G \to X.
\end{equation}
When $\X=X$ is a manifold (i.e.\ $G$ acts on $M$ locally freely), this map is a fibration with fiber $BG$.  The map $q$ induces
a map in cohomology, 
\begin{equation}\label{eq:quotientmap}
q^*: H^*(X;R) \longrightarrow H^*((M\times EG)/G;R) = H^*_G(M;R).
\end{equation}
This induced map is an isomorphism when
\begin{enumerate}
\item[1.] $G$ acts freely on $M$;
\item[2.] $G$ is a finite group and $R$ is a ring in which $|G|$ is invertible; and
\item[3.] $G$ acts locally freely on $M$ and $R$ is a field of characteristic $0$.
\end{enumerate}
This last item implies that the cohomology of the orbifold differs from the 
cohomology of its coarse moduli only in its torsion.  Notably,   when $R=\Z$, this ring does
in fact distinguish between the orbisphere in Figure~\ref{fig:orbisphere} and
a smooth sphere.

\begin{figure}[h]
  \begin{center}
  \psfrag{A}{$\cdots$}
    \epsfig{file=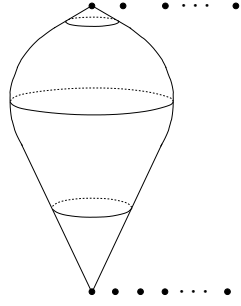,height=4cm}
    \caption{The inertia orbifold for the the orbisphere with two
  orbifold points as in Figure~\ref{fig:orbisphere}.  Each of the $p-1$
  points to the right of the north pole represents a $[
  \mathrm{pt}/\Z_p]$ and each of the $q-1$ to the right of the south
  pole a $[ \mathrm{pt}/\Z_q]$.}\label{fig:inertia} 
  \end{center}
 \end{figure}

The third invariant was introduced by Chen and Ruan \cite{CR:orbH} to
explain mathematically the {\bf stringy Betti numbers} and {\bf
stringy Hodge numbers} that physicists have attached to orbifolds.
To define this third invariant, we need to introduce the {\bf first
inertia orbifold} 
\begin{equation}
\label{ }
I^1(\X) := \left\{ \left(x,(g)_{\Gamma_x}\right)\ \Big| \ x\in \X \mbox{
and } (g)_{\Gamma_x} \mbox{ is a conjugacy class in }
\Gamma_x\right\}.
\end{equation}
This is again an orbifold, and $\X$ is the suborbifold called the {\bf identity sector} whose pairs consist of 
a point in $\X$ together with the identity element coset.   The other connected components of $I^1(\X)$ are called the {\bf twisted sectors}.   For a global quotient $\X=[M/G]$ with $G$
abelian, we may identify $I^1(\X) = \coprod_{g\in G} [M^g/G]$.  On the other hand, when $\X=[M/G]$ is  global quotient with $G$ finite, we may identify $I^1(\X) = \coprod_{g\in T} [M^g/C(g)]$, where the union is over $T$ a set of  representatives of conjugacy classes in $G$.  For
the orbisphere example, the inertia orbifold is shown in
Figure~\ref{fig:inertia}. 
 
\begin{definition}[\cite{CR:orbH}]
Given a presentation of an orbifold $\X = [M/G]$ as a global quotient,
the {\bf Chen-Ruan orbifold cohomology} of $\X$, as a vector space, is
defined to be the cohomology of the first inertia orbifold, 
\begin{equation}
\label{ }
H_{CR}(\X;R) := H(I^1(\X);R),
\end{equation}
with coefficients in a commutative ring $R$.  
\end{definition}

The Chen-Ruan ring is endowed with a $\Q$ grading, different from the
grading coming from singular cohomology.  For a connected component $\ZZ$ of
$I^1(\X)$ that lies in the $(g)$ piece of the first inertia, the
grading of $H(\ZZ;R)$ is shifted by a rational number which is twice
the {\bf age} of $\ZZ$.  The age is determined by the weights of the
action of the group element $g$ on the normal bundle to $\ZZ$ inside of
$\X$.  This is precisely where we need $\X$ to be (stably) almost
complex.  These  rational shifts ensure that the ranks of the
Chen-Ruan cohomology groups agree with the stringy invariants of an
orbifold. The dependence of the rational shifts on the normal bundles 
$\nu(\ZZ\subset\X)$ means that the Chen-Ruan ring is not in general functorial
for arbitrary morphisms of orbifolds.

To define a product on the Chen-Ruan cohomology, we must define higher
inertia.  The $\mathbf{n^{\mathrm{th}}}$
{\bf inertia orbifold} consists of tuples, a point in the orbifold and an $n$-tuple
of conjugacy classes. Restricting to the $2^{\mathrm{nd}}$ inertia,
there are natural maps 
\begin{equation}
\label{ }
e_1,e_2,\overline{e}_3 : I^2(X)\longrightarrow I^1(X)
\end{equation}
defined by 
\begin{eqnarray}
e_1(p,(g)_{\Gamma_p},(h)_{\Gamma_p}) & =&  (p, (g)_{\Gamma_p}),\\
e_2(p,(g)_{\Gamma_p},(h)_{\Gamma_p}) & =&  (p,(h)_{\Gamma_p}),\mbox{
and}\\
\overline{e}_3(p,(g)_{\Gamma_p},(h)_{\Gamma_p}) & =&  (p,
((gh)^{-1})_{\Gamma_p})
\end{eqnarray}
Chen and Ruan define the product of two classes $\alpha,\beta\in
H_{CR}(\X;R)$ to be
be
\begin{equation}
\label{ }
\alpha\smile\beta := (\overline{e}_3)_*(e_1^*\alpha\cup
e_2^*\beta\cup\varepsilon),
\end{equation}
where $e_1^*$ and $e_2^*$ are pull-back maps, $(\overline{e}_3)_*$ is
the push-forward, $\cup$ is the usual cup product, and $\varepsilon$
is the Euler class of the {\bf obstruction bundle}.  This Euler class
should be viewed as a quantum correction term.  It ensures that the
product respects the $\Q$-grading and that the product is
associative.  Neither of these properties is immediately obvious, and
proving the latter requires a rather substantial argument.  Since we can
avoid mention of the obstruction bundle in our computations, we
suppress further details here, and refer the curious reader to
\cite{AGV:orbifoldquantumproducts, CR:orbH, GHK:preorb} for additional
information.

The Chen-Ruan cohomology ring is the degree $0$ part of the (small) quantum cohomology ring.  
Hence, it is generally the most difficult of these
three invariants to compute.  It has been computed for orbifolds that
are global quotients of a manifold by a finite group in
\cite{FG:globalquotients}.  The definition was extended to the
algebraic category in \cite{AGV:orbifoldquantumproducts}, and in
\cite{BCS:toricvarieties}, the ring is computed for toric
Deligne-Mumford stacks with $\Q$, $\R$ and $\C$ coefficients. In the next
section, we review how to compute this ring for abelian symplectic
quotients, as demonstrated in \cite{GHK:preorb}.  In the last section, we will 
compute each of these invariants explicitly for weighted projective spaces.

We conclude this section with a table of these rings for an
orbisphere that has a $\Z_2$ singularity at the north pole and is
otherwise smooth.  This is an example of a {\bf weighted projective
space}, and is denoted $\C P^1_{1,2}$. 
\begin{equation}
\label{ }
\begin{array}{|c|c|c|}\hline
\X = [\C P^1_{1,2}] & \mbox{Ring} & \mbox{Grading} \\ \hline
H^*(X;\Z) & \Z[x]/\langle x^2\rangle & \deg(x)=2 \\ \hline
H^*(\X;\Z) & \Z[x]/\langle 2x^2\rangle & \deg(x)=2 \\ \hline
H_{CR}^{\phantom{i}*\phantom{i}}(\X;\Z) &  \Z[x,u]/\langle 2x^2, 2xu, 
u^2-x \rangle & \deg(x)=2,\ \deg(u)=1 \\ \hline
\end{array}
\end{equation}

\section{Why the symplectic category is convenient}\label{sec:symp}

In the past thirty years, tremendous progress has been made in
understanding the equivariant topology of Hamiltonian $T$-spaces and
its relationship to the ordinary topology of their quotients.  
For a compact torus $T=(S^1)^n$, the classifying bundle $ET$ is an
$n$-fold product of infinite dimensional spheres, and the classifying
space $BT$ is an $n$-fold product of copies of $\C P^\infty$.  Thus,
\begin{equation}
\label{ }
H_T^*(pt;\Z)= H^*(BT;\Z)=\Z[x_1,\dots,x_n],
\end{equation}
where $\deg(x_i)=2$.  The key
ingredient to understanding topology of Hamiltonian $T$-spaces is the
moment map.  Frankel  \cite{frankel} proved that for a Hamiltonian
$T$-action on a K\"ahler manifold, each component $\phi^\xi$ of the
moment map is a Morse-Bott function on $M$, and generically the
critical set is the fixed point set $M^T$.  In his paper
\cite{At:convexity} on the Convexity Theorem~\ref{thm:convex}, Atiyah
generalized this work to the purely symplectic setting.

Building on the work of Frankel and Atiyah,  Kirwan developed techniques to
prove two fundamental theorems
that allow us to understand the cohomology of Hamiltonian $T$-spaces
and their quotients.  The first is a version of {\bf localization}: it
allows us to make global computations by understanding fixed point
data. While this theorem is not explicitly stated in her book \cite{K:quotients}, it 
does follow immediately from her work in Chapter~5.

\begin{thmi}[\cite{K:quotients}]\label{th:inj}
Let $M$ be a compact Hamiltonian $T$-space.
The inclusion map $M^T\hookrightarrow M$ induces
\begin{equation}
\label{ }
i^*:H_T^*(M;\Q) \longrightarrow H_T^*(M^T;\Q)
\end{equation}
an {\bf injection} in equivariant cohomology.
\end{thmi}

The compactness hypothesis is stronger than strictly necessary.  We
may replace it with a properness condition on the moment map.  The
proof relies on the fact that a generic component of the moment
map is an {\bf equivariantly perfect} Morse-Bott function on $M$. The
image of this injection has been computed in many examples, including
toric varieties, coadjoint orbits of compact connected semisimple Lie
groups, and coadjoint orbits of Kac-Moody groups.  These computations
initially appeared in \cite{CS:injectivity,GKM} and further
generalizations are described in \cite{GoH,GuH,GZ,HHH}. 

The second theorem relates the equivariant topology of a Hamiltonian
$T$-space to the ordinary topology of its reduction.

\begin{thms}[\cite{K:quotients}]\label{th:surj}
Let $M$ be a compact Hamiltonian $T$-space, and $\alpha$ a regular 
value of the moment map.
The inclusion $\Phi^{-1}(\alpha)\hookrightarrow M$ induces
\begin{equation}
\label{ }
\kappa: H_T^*(M;\Q) \longrightarrow H_T^*(\Phi^{-1}(\alpha);\Q)
\end{equation}
a {\bf surjection} in equivariant cohomology.
\end{thms}

Again for surjectivity, compactness is more than is necessary.  Most
importantly, this result does apply to linear actions of a torus on
$\C^d$ with a proper moment map.   
The key idea in the proof is to use the function $||\Phi-\alpha| |^2$ as a {\em Morse-like} function, now known as a {\bf Morse-Kirwan function}.  The critical sets are not non-degenerate, but one may still explicitly understand them via a local normal form. It is then possible to prove that $||\Phi-\alpha| |^2$ is an equivariantly perfect function on $M$.

The kernel of the map $\kappa$
can be computed using methods in
\cite{Go:effective,JK:kernel, TW:symplecticquotients}.  Using the fact that at a
regular value,
$H_T^*(\Phi^{-1}(\alpha);\Q)\cong H^*(M/\!\!/T(\alpha);\Q)$, we have
a diagram 
\begin{equation}
\label{ }
\begin{array}{c}
\xymatrix{
0 \ar[r] & \ker(\kappa)\ar@{^{(}->}[r] & H_T^*(M;\Q)\ar@{->>}[r]^(0.4){\kappa}  )\ar@{^{(}->}[d]^{i^*} & H^*(M/\!\!/T(\alpha);\Q) \ar[r] & 0\\
& & H_T^*(M^T;\Q) & & 
}
\end{array}.
\end{equation}
Thus, by computing $\im(i^*)$ and $\ker(\kappa)$, we may derive an
explicit presentation of the cohomology $H^*([M/\!\!/T(\alpha)];\Q)$
of an orbifold arising as a symplectic quotient.

We now turn to a generalization of Theorem~\ref{th:surj} in the context of orbifolds and the Chen-Ruan ring. 

\begin{thmcr}[\cite{GHK:preorb}]\label{th:main}
Let $M$ be a compact Hamiltonian $T$-space, and $\alpha$ a regular value of the moment map.
The inclusion $\Phi^{-1}(\alpha)\hookrightarrow M$ induces
\begin{equation}\label{eq:cr-surj}
\mathcal{K}: \bigoplus_{g\in T} H_T^*(M^g;\Q) \longrightarrow \bigoplus_{g\in T} H_T^*(\Phi^{-1}(\alpha)^g;\Q)
\end{equation}
a {\bf surjection}.  Moreover, these are $\R\times T$-graded rings, $\mathcal{K}$ is a map of graded rings, and there is an isomorphism of graded rings
\begin{equation}\label{eq:CRiso}
\bigoplus_{g\in T} H_T^*(\Phi^{-1}(\alpha)^g;\Q) \cong H_{CR}^{\phantom{i}*\phantom{i}}(M/\!\! /T(\alpha);\Q).
\end{equation}
\end{thmcr}

The surjectivity \eqref{eq:cr-surj} is a direct consequence of the Surjectivity
Theorem~\ref{th:surj} applied to each space $M^g$ (and again the compactness is not strictly
necessary).  The hard work is defining the grading and ring structure, and proving that $\mathcal{K}$ is a map of graded rings.  Goldin, Knutson, and the author define a ring structure $\smile$
that generalizes a definition (for $G$ finite) of Fantechi and G\"ottsche \cite{FG:globalquotients}.  Using
this definition, we may deduce \eqref{eq:CRiso}; on the other hand, associativity of
this product is not at all obvious.
Making use of the injection $i^*$ on each piece, there is
an alternative 
product $\star$ on the ring $\bigoplus_{g\in T}
H_T^*(M^g;\Q)$ that is much simpler to compute, and clearly
associative.   This alternative product has the advantage that 
it avoids all mention of the obstruction bundle; instead it relies only
on fixed point data (i.e.\ the topology of the fixed point set and isotropy data for the action of the torus on the normal bundles to the fixed point components).

Another key point is that although the ring on the left of \eqref{eq:cr-surj} is quite
large, there is a finite subgroup $\Gamma$ of $T$, generated by all
finite-order elements that stabilize some regular point in $M$, so that the
${\mathbf \Gamma}${\bf-subring} 
\begin{equation}\label{eq:gammapiece}
\bigoplus_{g\in \Gamma} H_T^*(M^g;\Q)
\end{equation}
still surjects
onto the Chen-Ruan cohomology of the reduction.  Thus, while it
appears that we have made the computation much more complicated, it
turns that there is still an effective algorithm to complete it. 
For full details, please refer to \cite{GHK:preorb}.  
We now return to our examples.

\begin{example}\label{EG:symptoricorb}
For a {\bf symplectic toric orbifold}  $\C^n/\!\!/K (\alpha)$,
we may use the combinatorics of its labeled polytope to
establish an explicit presentation of the Chen-Ruan cohomology of
these orbifolds  \cite[\S\ 9]{GHK:preorb}. In the cases where the symplectic picture is
identical to the algebraic, the \cite{GHK:preorb} results replicate
those of \cite{BCS:toricvarieties}.  
\end{example}

\begin{example}
For the $T$-action on a coadjoint orbit $M=
\mathcal{O}_\lambda$, the symplectic reduction $M/\!\!/T(\alpha)$ is a
{\bf weight variety}.  The equivariant cohomology of $M$ may be read
directly from its moment polytope, as may the orbifold singularities of the
reduction.  In this case, the Theorem~\ref{th:main} yields an explicit
combinatorial description of the Chen-Ruan cohomology of the weight
variety. 
\end{example}

We conclude this section with a brief remark on coefficients.  In both
Theorems~\ref{th:inj} and \ref{th:surj}, the rational coefficients are
necessary.  For the Injectivity Theorem~\ref{th:inj}, we may prove the result over $\Z$ with
the additional hypothesis that $H^*(M^T;\Z)$ contains no torsion.
The Surjectivity Theorem~\ref{th:surj} over $\Z$ requires much stronger hypotheses.  We will see
in the next section that we may compute integrally for weighted
projective spaces, but that a simple product of two weighted
projective spaces yields a counter-example to the general theorem.
Tolman and Weitsman verify surjectivity over $\Z$ for a rather
restrictive class of torus actions \cite{TW:symplecticquotients}.
This topic is being more closely examined for a larger collection
of actions by Susan Tolman and the author \cite{HT}.

\section{The case of weighted projective spaces}\label{sec:wps}

Let $b = (b_0,\dots,b_n)$ be an $(n+1)$-tuple of positive integers.
Consider the circle action on $\C^{n+1}$ given by 
\begin{equation}
\label{ }
t\cdot (z_0,\dots,z_n) = (t^{b_0}z_0,\dots,t^{b_n}z_n)
\end{equation}
for each $t\in S^1$.  This action preserves the unit sphere $S^{2n+1}$,
and the {\bf weighted projective space} $\C P^n_{(b)}$ is the quotient
of $S^{2n+1}$ by this locally free circle action.  This is a
symplectic reduction because $S^{2n+1}$ is (up to equivariant
homeomorphism) a regular level set for a moment map $\Phi$ for the
weighted $S^1$ action on $\C^{n+1}$.   Thus, 
\begin{equation}
\label{ }
\C P^n_{(b)} \cong \C^{n+1}/\!\! /S^1.
\end{equation}
Nonetheless, we may continue our analysis without invoking the full symplectic machinery: the
arguments simplify greatly in this special case.

When the $b_i$ are relatively prime (that is, $\gcd(b_0,\dots,b_n)=1$), this fits into the
framework of symplectic toric orbifolds discussed in
Example~\ref{EG:symptoricorb}.  When the $b_i$ are not relatively
prime, we let $g = \gcd(b_0,\dots,b_n)$, and note that there is a global $\Z_g$ stabilizer. 
In this case, the orbifold $[\C P^n_{(b)}]$ is not
reduced.  Its coarse moduli space is the same as the coarse moduli space of $\C P^n_{(b/g)}$, where
$(b/g)$ denotes the sequence of integers $(b_0/g,\dots,b_n/g)$; and as a non-reduced orbifold, it corresponds to a {\bf gerbe} 
$$
[\C P^n_{(b)}]\to[\C P^n_{(b/g)}].
$$   
It is important to include the case when the
$b_i$ are not relatively prime, because such non-reduced weighted projective spaces may well show up as suborbifolds of a reduced weighted projective space.

\medskip

\noindent {\bf The cohomology of the topological space $\C P^n_{(b)}$.}  Kawasaki
studied the singular cohomology ring of (the coarse moduli space of) weighted projective spaces 
\cite{K:weightedproj}.  To present the product structure, we will need the integers
\begin{equation}
\label{ }
\ell_k = \ell_k^{(b)} := \lcm\left\{ \frac{b_{i_0}\cdots
b_{i_k}}{\gcd(b_{i_0},\dots,b_{i_k})}\ \bigg| \ 0\leq i_0<\cdots
<i_k\leq n \right\},
\end{equation}
for each $1\leq k\leq n$.

\begin{theorem}[\cite{K:weightedproj}]\label{th:kaw}
The integral cohomology of $\C P^n_{(b)}$ is 
\begin{equation}
\label{ }
H^i(\C P^n_{(b)};\Z) \cong \left\{\begin{array}{ll}
	\Z & \mbox{if } i=2k,\ 0\leq k\leq n, \\
	0 & \mbox{otherwise.}
\end{array}\right.
\end{equation}
Moreover, letting $\gamma_{i}$ denote the generator of $H^{2i}( \C P^n_{(b)};\Z)$,
we have
\begin{equation}
\label{ }
\gamma_k\cup\gamma_m = \frac{\ell_k\cdot\ell_m}{\ell_{m+k}} \gamma_{m+k}.
\end{equation}
\end{theorem}

\noindent {\textsc{Outline of the Proof}.}
Let $G_k$ denote
the group of $k^{\mathrm{th}}$ roots of unity.  Then as a topological
space, $\C P^n_{(b)}$ is
homeomorphic to a quotient of ordinary projective space $\C P^n$ by
the finite group
\begin{equation}
\label{ }
G_{(b)} = G_{b_0}\times\cdots\times G_{b_n}.
\end{equation}
Explicitly, using standard homogeneous coordinates, the map
\begin{eqnarray}
p_b: \C P^n & \longrightarrow & \C P^n_{(b)} \\
\left[ z_0:\cdots:z_n\right] & \longmapsto & [z_0^{b_0}:\cdots: z_n^{b_n}]
\end{eqnarray}
induces the homeomorphism $\C P^n/G_{(b)}\cong\C P^n_{(b)}$.
This then induces an isomorphism in
singular cohomology with rational coefficients, since over $\Q$ we have the
isomorphisms
\begin{equation}
\label{ }
H^*(\C P^n;\Q) \cong H^*(\C P^n;\Q)^{G_{(b)}} \cong H^*(\C
P^n/G_{(b)};\Q). 
\end{equation}

Over the integers, the computation is a bit more subtle.  Using twisted lens spaces, 
Kawasaki verifies that  just as for $\C P^n$, the 
cohomology ring $H^*(\C P^n_{(b)};\Z)$ is torsion-free with a copy of $\Z$ in 
each even degree between $0$ and $2n$;
however the product structure is twisted by the weights $b_i$.  Moreover, there are cases
when we need all $n$ generators $\gamma_1,\dots,\gamma_n$ to present this ring.

Kawasaki  showed that the map
\begin{equation}
\label{ }
p_b^*: H^{2k}(\C P^n_{(b)};\Z)\longrightarrow H^{2k}(\C P^n;\Z)
\end{equation}
is multiplication by $\ell_k$ for all $1\leq k\leq n$.  From this, we may deduce that
\begin{equation}
\gamma_1\cup\gamma_k = \frac{\ell_1\cdot \ell_{k+1}}{\ell_k}
\gamma_{k+1}, \mbox{ and hence }\  \gamma_k\cup\gamma_m = \frac{\ell_k\cdot\ell_m}{\ell_{m+k}} \gamma_{m+k}.
\end{equation}
The result now follows.  \hfill\qed

\medskip

It is important to note that $p_b$ is {\bf not} an isomorphism of orbifolds.  Indeed, the isotropy group 
at any point in $\C P^n_{(b)}$ is the stabilizer group of any lift of the point in $S^{2n+1}$.  Thus,
all isotropy groups for $\C P^n_{(b)}$ are cyclic.  On the other hand, the orbifold 
$\C P^n/G_{(b)}$ has points with 
isotropy group $G_{(b)}$, which may not be cyclic.  Nevertheless 
 we will make use of the map $p_b$ to understand the structure of the cohomology
of the orbifold $[\C P^n_{(b)}] = [S^{2n+1}/S^1_{(b)}]$.  Here, the subscript $(b)$  on $S^1$ indicates that 
the circle action is  weighted by the integers $(b) = (b_0,\dots,b_n)$.

\medskip

\noindent {\bf The cohomology of the orbifold $[\C P^n_{(b)}]$.} 
Since the weighted projective space $[\C P^n_{(b)}]$ is a symplectic reduction, it is possible invoke the results from Section~\ref{sec:symp} to determine the cohomology of the orbifold $H^*([\C P^n_{(b)}];\Q)$, and to apply results from \cite[\S 9]{GHK:preorb} to obtain a presentation over $\Z$.  We give a direct argument here that is similar in spirit, but that avoids much of this big machinery; we then compare this to Kawasaki's Theorem~\ref{th:kaw}.

\begin{theorem}\label{th:surjZwps}
The cohomology of the orbifold $[\C P^n_{(b)}]$ is
\begin{equation}\label{eq:coh_orb}
H^*([\C P^n_{(b)}];\Z) = H_{S^1_{(b)}}^*(S^{2n+1};\Z) \cong \frac{\Z[u]}{\langle b_0\cdots b_n u^{n+1}\rangle}.
\end{equation}
Moreover, the natural map
\begin{equation}\label{eq:kaw-rel}
q_{(b)}^*:H^*(S^{2n+1}/S^1_{(b)};\Z)\longrightarrow H_{S^1_{(b)}}^*(S^{2n+1};\Z)
\end{equation}
is completely determined by $q_{(b)}^*(\gamma_1) = \ell_1 \cdot u = \lcm(b_0,\dots,b_n)\cdot u$.
\end{theorem}

\noindent {\textsc{Proof}.}
Consider the (weighted) circle action of $S^1$ on $\C^{n+1}$ given by
\begin{equation}
\label{ }
t\cdot (z_0,\cdots,z_n) = (t^{b_0}\cdot z_0,\cdots,t^{b_n}\cdot z_n).
\end{equation}
The unit sphere $S^{2n+1}$ is invariant under this action, so we get a long
exact sequence in $S^1$-equivariant cohomology for the pair $(\C^{n+1},S^{2n+1})$,
\begin{equation}
\label{ }
\cdots\to  H^i_{S^1_{(b)}}(\C^{n+1},S^{2n+1} ;\Z)  \stackrel{\alpha}{\to} H^i_{S^1_{(b)}}(\C^{n+1} ;\Z) \stackrel{\beta}{\to} H^i_{S^1_{(b)}}(S^{2n+1} ;\Z) \to \cdots.
\end{equation}
Thinking of $(\C^{n+1},S^{2n+1})$ as a disk and sphere bundle over a point, we may use the Thom isomorphism to identify $H^i_{S^1_{(b)}}(\C^{n+1},S^{2n+1} ;\Z)  \cong H^{i-2(n+1)}_{S^1_{(b)}}(\C^{n+1} ;\Z)$. Under this identification, the map $\alpha$ is the cup product with the equivariant Euler class
\begin{equation}
\label{ }
e_{S^1_{(b)}}(\C^{n+1}) = b_o\cdots b_n u^{n+1}.
\end{equation}
Thus, the map $\alpha$ is injective, so the long exact sequence splits into short exact sequences
\begin{equation}
\label{eq:ses}
0\to H^i_{S^1_{(b)}}(\C^{n+1},S^{2n+1} ;\Z)  \stackrel{\alpha}{\to} H^i_{S^1_{(b)}}(\C^{n+1} ;\Z) \stackrel{\beta}{\to} H^i_{S^1_{(b)}}(S^{2n+1} ;\Z) \to 0.
\end{equation}
Thus, we have a surjection
\begin{equation}
\label{ }
\beta: \Z[u]\cong H^*_{S^1_{(b)}}(\C^{n+1} ;\Z) \longrightarrow H^*_{S^1_{(b)}}(S^{2n+1} ;\Z). 
\end{equation}
Moreover, the exactness of \eqref{eq:ses} means that the kernel of $\beta$ is equal to the image of $\alpha$, namely all multiples of the equivariant Euler class.  This establishes \eqref{eq:coh_orb}.

Turning to \eqref{eq:kaw-rel}, the map 
\begin{equation}
q_{(b)}^*: H^*(S^{2n+1}/S^1_{(b)};\Z)\longrightarrow H_{S^1_{(b)}}^*(S^{2n+1};\Z),
\end{equation}
is exactly the one defined in \eqref{eq:quotientmap}.  We know that this is an isomorphism over $\Q$, so $q_{(b)}^*$ must map $\gamma_1$ to a multiple of $u$.  Moreover, because $b_0\cdots b_nu^{n+1}$ is zero in $H_{S^1_{(b)}}^*(S^{2n+1};\Z)$, we must have that 
\begin{equation}
\label{ }
(q_{(b)}^*(\gamma_1))^{n+1}\in \left\langle b_0\cdots b_n u^{n+1}\right\rangle.
\end{equation}

To determine the image of the class $\gamma_1$,
we return to the map $p_b : \C P^n \to \C P^n_{(b)}$.  This map lifts to  maps 
on $S^{2n+1}$ and $\C^{n+1}$ given by
\begin{equation}
\label{ }
\begin{array}{c}
\xymatrix{
\C^{n+1}  \ar[r]^{\Pi_b} & \C^{n+1} \\
 S^{2n+1}\ar[r]^{\pi_b}\ar@{^(->}[u]^{i} & S^{2n+1} \ar@{_(->}[u]_{i} 
 }
\end{array}.
\end{equation}
\begin{eqnarray}
\phantom{BOoooO} 0\neq( z_0,\cdots,z_n) & \longmapsto  & \frac{\sum |z_i|^2}{\sum |z_i^{b_i}|^2}\cdot \left(z_0^{b_0},\cdots, z_n^{b_n}\right) \\
(0,\dots,0) & \longmapsto & (0,\dots,0)
\end{eqnarray}
The maps in this diagram are all equivariant with respect to the standard circle action on the
left-hand spaces and the $(b)$-weighted circle action on the right-hand spaces.  Thus, we have 
a diagram of maps
\begin{equation}
\label{ }
\begin{array}{c}
\xymatrix{
(\C^{n+1}\times S^1)/S^1  \ar[r]^{\Pi_b} & (\C^{n+1} \times S^1_{(b)})/S^1_{(b)} \\
( S^{2n+1}\times S^1)/S^1 \ar[r]^{\pi_b}\ar@{^(->}[u]^{i}\ar[d]^{q} & (S^{2n+1} \times S^1_{(b)})/S^1_{(b)} \ar@{_(->}[u]_{i} \ar[d]_{q_{(b)}} \\
\C P^n\ar[r]_{p_{b}} &  \C P^n_{(b)}
} 
\end{array}.
\end{equation}
Applying singular cohomology $H^*(\phantom{-};\Z)$ and identifying equivariant cohomology, we have a commutative diagram
\begin{equation}
\label{ }
\begin{small}
\begin{array}{c}
\xymatrix{
\Z[x] \ar@{=}[r] & H_{S^1}^*(\C^{n+1};\Z)  \ar@{<-}[r]^{\Pi^*_b} & H_{S^1_{(b)}}^*(\C^{n+1};\Z) & \Z[u]\ar@{=}[l] \\
\frac{\Z[x]}{\left\langle x^{n+1}\right\rangle} \ar@{=}[d]\ar@{=}[r] & H_{S^1}^*(S^{2n+1};\Z) \ar@{<-}[r]^{\pi_b^*}\ar@{<-}[u]^{\kappa}\ar@{<-}[d]^{q^*} & H_{S^1_{(b)}}^*(S^{2n+1};\Z) \ar@{<-}[u]_{\kappa_{(b)}} \ar@{<-}[d]_{q^*_{(b)}} & \frac{ \Z[u]}{\left\langle b_0\cdots b_nu^{n+1}\right\rangle} \ar@{=}[l] \\
\frac{\Z[x]}{\left\langle x^{n+1}\right\rangle} \ar@{=}[r] & H^*(\C P^n;\Z) \ar@{<-}[r]_{p^*_{b}} &  H^*(\C P^n_{(b)};\Z) & \Z\oplus\Z\gamma_1\oplus\cdots\oplus\Z\gamma_n \ar@{=}[l]
} 
\end{array}.
\end{small}
\end{equation}
 Because $\C^{n+1}$ equivariantly deformation retracts to a point, the map $\Pi_b^*$ maps the 
generator $u$ to $x$.  The commutativity of the top square then implies that $\pi_b^*(u) = x$.  Thus, we know that
\begin{eqnarray}
\pi_b^*(q^*_{(b)}(\gamma_1) & = & q^*(p_b^*(\gamma_1))\\
& = & q^*(\ell_1 \cdot x), \ \ \  \mbox{ by Kawasaki's result,}\\
& = & \ell_1\cdot x, \ \ \ \ \ \ \ \ \mbox{ since } q^*  \mbox{ is an equality,} \\
& = & \pi_b^*( \ell_1\cdot u).
\end{eqnarray}
In low degree, $\pi_b^*$ is injective, so we may conclude that $q_{(b)}^*(\gamma_1) = \ell_1 \cdot u$. Noting that $\ell_1=\lcm(b_0,\dots,b_n)$ completes the proof.\hfill\qed

\medskip

Over the integers, this invariant does distinguish a weighted projective space from the standard one; however, it may not differentiate between two weighted projective spaces.  For example, the cohomology rings of the orbifolds $[\C P^1_{2,2}]$ and $[\C P^1_{4,1}]$ are identical.  They are both
\begin{equation}
\label{ }
\frac{\Z[u]}{\langle 4u^2\rangle}.
\end{equation}

We note that these surjectivity techniques do not generally work over the integers.  To see this, we note that for any abelian reduction of affine space, the domain of the Kirwan map $H_T^*(\C^N;\Z)$ has terms only in even degrees.  If we consider the simple product $[\C P^1_{1,2}\times \C P^1_{1,2}]$, we may compute the cohomology of this orbifold using the above result and the K\"unneth formula.  Since $[\C P^1_{1,2}]$ has $2$-torsion in high degrees, the $\mathrm{Tor}$ term from the K\"unneth formula plays a role, yielding $2$-torsion in high odd degrees in the cohomology of the orbifold $[\C P^1_{1,2}\times \C P^2_{1,2}]$.  Thus, surjectivity must fail over the integers in this example.  We note that any failure over $\Z$ must be due to problems with torsion, because surjectivity does hold over $\Q$.

\medskip

\noindent {\bf The Chen-Ruan orbifold cohomology of $[\C P^n_{(b)}]$.} When computing the Chen-Ruan ring, it is important to recall that a weighted projective space is a circle reduction.  Thus, the finite group $\Gamma$ for which the $\Gamma$-piece surjects onto  $H_{CR}^{\phantom{i}*\phantom{i}}([\C P^n_{(b)}];\Z)$ is a cyclic group.  For any vector $v$ that is non-zero is a single coordinate, say the $i^{\mathrm{th}}$ coordinate, the stabilizer of $v$ is $\Z_{b_i}$.  Thus, the group $\Gamma$ generated by all finite stabilizers is the $\ell^{\mathrm{th}}$ roots of unity $\Z_\ell \subset S^1$, where $\ell=\lcm(b_0,\dots,b_n)$.  Hence, we have a surjection
\begin{equation}
\label{ }
\xymatrix{
\Z[u,\alpha_0,\alpha_1,\dots,\alpha_{\ell-1}] \ar@{->>}[r] & H_{CR}^{\phantom{i}*\phantom{i}}([\C P^n_{(b)}];\Z).
}
\end{equation}
In this case, thinking of $e^{\frac{2\pi i k}{\ell}} = \zeta_k\in \Z_\ell \subset S^1$, $\alpha_k$ denotes a generator for  
\begin{equation}
\label{ }
H^*_{S^1_{(b)}}((\C^{n+1})^{\zeta_k};\Z).
\end{equation}
To complete the computation, we must determine the orbifold product 
\begin{equation}
\label{ }
\alpha_i\smile\alpha_j=\alpha_i\star\alpha_j
\end{equation}
and the kernel of the orbifold Kirwan map \eqref{eq:cr-surj}.  For any integer $m\in \Z$, we let $[m]$ denote the smallest non-negative integer congruent to $m$ modulo $\ell$.  For any rational number $q\in\Q$, $\langle q\rangle_f$ denotes its fractional part.  Finally, we let 
\begin{equation}
\label{ }
a_k(m) := \frac{[ b_k\cdot m]}{\ell} = \left\langle  \frac{b_k\cdot m}{\ell} \right\rangle_f.
\end{equation}
This is the rational number such that $\zeta_m$ acts on the $k^{\mathrm{th}}$ coordinate by $e^{2\pi i a_k(m)}$.

\begin{theorem}\label{th:cr}
The Chen-Ruan orbifold cohomology of $\C P^n_{(b)}$ is
\begin{equation}
\label{ }
H_{CR}^{\phantom{i}*\phantom{i}}([\C P^n_{(b)}];\Z) \cong \frac{\Z[u,\alpha_0,\alpha_1,\dots,\alpha_{\ell-1}] }{\mathcal{I}+\mathcal{J}},
\end{equation}
where $u$ is a class in degree $2$, 
\begin{equation}
\label{ }
\deg(\alpha_j) = 2 \sum_{k=0}^n a_k(j).
\end{equation}
Here, $\mathcal{I}$ is the ideal
\begin{equation}\label{eq:product}
\mathcal{I} = \left\langle \alpha_i\alpha_j-\left(\prod_{k=0}^n (b_ku)^{a_k(i)+a_k(j)-a_k(i+j)}\right)\alpha_{[i+j]}\right\rangle
\end{equation}
generated by the $\star$ product structure, and $\mathcal{J}$ is 
\begin{equation}
\label{ }
\mathcal{J} = \sum_{j=0}^n \left\langle\left( \prod_{a_k(j)=0} b_ku\right) \alpha_i\right\rangle,
\end{equation}
the kernel of the surjection $\mathcal{K}$ of the orbifold Kirwan map.
\end{theorem}

\begin{remark}
The generator $u$ is the generator of $S^1$-equivariant cohomology and hence has degree $2$.  The generator $\alpha_k$ is a placeholder for the cohomology of the $\zeta_k$-sector.
\end{remark}

\begin{remark}
Note that the generator $\alpha_0$ is the placeholder for the identity sector.  Indeed, we always have
\begin{equation}
\label{ }
\alpha_0\star\alpha_0 =\alpha_0
\end{equation}
as a consequence of the relation in \eqref{eq:product} where $i=j=0$, hence  we may think of $\alpha_0$ as $1$.
\end{remark}

\begin{remark}
The reader may use this theorem to check that $H_{CR}^{\phantom{i}*\phantom{i}}(\phantom{-};\Z)$ does distinguish $[\C P^1_{2,2}]$ from $[\C P^1_{4,1}]$.
\end{remark}

\noindent {\textsc{Proof}.}
We use the $\star$ product given by Equation~(2.1) in \cite{GHK:preorb}.  In the case of a weighted circle action on $\C^{n+1}$, there is exactly one fixed point (the origin), any generator $\alpha_i$ restricted to that fixed point is $1$, and the equivariant Euler class for the $k^{\mathrm{th}}$ coordinate is precisely $b_ku$, whence
\begin{equation}
\label{ }
 \alpha_i\star\alpha_j=\left(\prod_{k=0}^n (b_ku)^{a_k(i)+a_k(j)-a_k(i+j)}\right)\alpha_{[i+j]}.
\end{equation}

Turning to the kernel computation, each $(\C^{n+1})^{\zeta_j}$ has a weighted $S^1$ action, and so we apply Theorem~\ref{th:surjZwps} to this subspace.  Thus, for the $\zeta_j$-sector,  the kernel contribution is the equivariant Euler class of $(\C^{n+1})^{\zeta_j}$ times the placeholder $\alpha_j$.  We note that $(\C^{n+1})^{\zeta_j}$ contains the $k^{\mathrm{th}}$ coordinate subspace precisely when $a_k(j)=0$.  Hence,
\begin{equation}
\label{ }
e_{S^1_{(b)}}((\C^{n+1})^{\zeta_j}) = \prod_{a_k(j)=0}^{\phantom{o}} b_ku,
\end{equation}
and the theorem follows.  \hfill\qed

\medskip

This theorem is an immediate consequence of Theorem~\ref{th:surjZwps} and \cite{GHK:preorb}.  The importance of this description is its ease in computation, since it avoids any computation of a labeled moment polytope ({\em \'a la} \cite{LT:toricorbifolds}) or of a stacky fan ({\em \'a la}  \cite{BCS:toricvarieties}).  We demonstrate this computational facility in the following concluding example.

\begin{example}
Consider the weighted projective space $[\C P^5_{1,2,2,3,3,3}]$.  This is a symplectic reduction of $\C^6$, the group $\Gamma$ is $\Z/6\Z$, and so the Chen-Ruan orbifold cohomology of $[\C P^5_{1,2,2,3,3,3}]$ is a quotient of
\begin{equation}
\label{ }
\Z[u,\alpha_0,\alpha_1,\alpha_2,\alpha_3,\alpha_4,\alpha_5].
\end{equation}
The following chart contains the data needed to compute the ideals $\mathcal{I}$ and $\mathcal{J}$.
\begin{equation}\label{eq:123table}
\begin{array}{c||c|c|c|c|c|c|}
g & \zeta_0 & \zeta_1 & \zeta_2 & \zeta_3 & \zeta_4 & \zeta_5 \\
\hline\hline (\C^6)^g & \C^6 & \{ 0\} & 3\C_{(3)} & 2\C_{(2)} & 3\C_{(3)} &
\{ 0\} \\ \hline \Year_{\C_{(1)}}(g) & 0 & \frac{1}{6} & \frac{1}{3} & \frac{1}{2} & \frac{2}{3} & \frac{5}{6}  \\
\hline \Year_{\C_{(2)}}(g) & 0 & \frac{1}{3} & \frac{2}{3}  & 0 &
\frac{1}{3} & \frac{2}{3} \\ \hline \Year_{\C_{(3)}}(g)  & 0 &
\frac{1}{2}  & 0 & \frac{1}{2}  & 0 & \frac{1}{2} \\ \hline 
2\cdot\age(g) & 0 & \frac{14}{3} & \frac{10}{3} & 4 & \frac{8}{3} & \frac{22}{3} \\ \hline
\genfrac{}{}{0pt}{0}{\mbox{generator of}}{H_{S^1_{(b)}}^*((\C^6)^g;\Z)} & \alpha_0 &
\alpha_1 & \alpha_2 &\alpha_3 & \alpha_4 & \alpha_5\\ \hline
e_{S^1_{(b)}}((\C^6)^g) & 108u^6& 1& 27u^3&4u^2 &27u^3 & 1\\ \hline
\end{array}
\end{equation}
Note that because of the multiplicities,
\begin{equation}
\label{ }
2\cdot \age(g) = 2\cdot \big[a_1(g) +2a_2(g)+3a_3(g)\big].
\end{equation}

Since $\alpha_0=1$, and $\alpha_1$ and $\alpha_5$ are in the kernel ideal $\mathcal{J}$, we only need to compute the products among $\alpha_2$, $\alpha_3$ and $\alpha_4$.  For example, we may compute
\begin{equation}
\label{ }
\alpha_2\star\alpha_2 = (u)^{\frac{1}{3}+\frac{1}{3}-\frac{2}{3}}\left( (2u)^{\frac{2}{3}+\frac{2}{3}-\frac{1}{3}}\right)^2 \left( (3u)^{0+0-0}\right)^3\alpha_4 = 4u^2\alpha_4.
\end{equation} 
All of the products  contributing to $\mathcal{I}$, then, are summarized in the following table.
\begin{equation}\label{eq:multtable}
\begin{array}{c||c|c|c|}
\star & \alpha_2 &  \alpha_3 &  \alpha_4 \\ \hline \hline
\alpha_2 & 4u^2\alpha_4 &  \alpha_5=0 & 4u^3 \\ \hline
 \alpha_3 & & 27u^4 & u\alpha_1=0 \\ \hline 
  \alpha_4 & & & u\alpha_2 \\ \hline
\end{array}.
\end{equation}
Thus, as a ring,
\begin{equation}
\label{ }
H_{CR}^{\phantom{i}*\phantom{i}}([\C P^n_{1,2,2,3,3,3}];\Z) \iso \frac{\Z[u,\alpha_0,\alpha_1,\alpha_2,\alpha_3,\alpha_4,\alpha_5]}{\mathcal{I}+\left\langle 108u^6, \alpha_1, 27u^3\alpha_2, 4u^2\alpha_3, 27u^3\alpha_4, \alpha_5  \right\rangle}.
\end{equation}
This generalizes Jiang's computation \cite{J:wps} to a computation over $\Z$.
\end{example}

\bibliographystyle{amsalpha}

\end{document}